\documentclass{article}
\usepackage{amssymb,amsmath,amsthm}

\def\={~=~}
\newcommand{\lp}{\left(}
\newcommand{\lb}{\left\lbrack}
\newcommand{\rp}{\right)}
\newcommand{\rb}{\right\rbrack}

\newtheorem{lemma}{Lemma}

\begin{document}

\title{A Class of Periodic Continued Radicals}
\author{Costas J. Efthimiou}
\date{}
\maketitle

\begin{abstract}
We compute the limits of a class of periodic continued radicals and we establish a connection between them and the fixed
points of the Chebycheff polynomials. 
\end{abstract}

\section{Introduction.}
Continued radicals 
{\footnotesize
$$
      a_0\sqrt{b_1+a_1\sqrt{b_2+a_2\sqrt{b_3+a_3\sqrt{b_4+\cdots}}}}    
$$
}
have been well known among mathematicians [1,3--6,9--11]
and they have 
even appeared  at mathematical
competitions. For example, Ramanujan's famous result \cite{Ramanujan}
{\footnotesize
$$
    \sqrt{1+2\sqrt{1+3\sqrt{1+4\sqrt{1+\cdots}}}}               \=  3
$$
}
 was given as one of the problems on the Putnam Mathematical Competition in 1966. However, the literature 
 on them is rather scant and, with only some exceptions,
 the main results  consider cases with positive numbers $a_i$ and $b_i$. Among the results on continued radicals with negative coefficients
 are  some  by Ramanujan himelf \cite{Ramanujan}, such as  the continued radical
 {\footnotesize
 $$
    \sqrt{a-\sqrt{a +\sqrt{a+\sqrt{a-\sqrt{a+\cdots}}}}}                
$$
 }
with period 3, and
Problem 1174 in  \cite{Dence} (motivated by a 1953 Putnam Problem for which $a=7$)
{\footnotesize
$$
    \sqrt{a-\sqrt{a +\sqrt{a-\sqrt{a+\sqrt{a-\cdots}}}}}                
$$
}
with period 2.

\section{The Problem.}
In this brief article we find the values for a class of periodic continued radicals of the form
{\footnotesize
\begin{equation}
         a_0\sqrt{2+a_1\sqrt{2+a_2\sqrt{2+a_3\sqrt{2+\cdots}}}}     ~,
\label{eq:OurRadical}
\end{equation}
}
where for some positive integer $n$,
$$
     a_{n+k} \= a_k~,~~~k=0,1,2,\dots~,
$$     
and
$$
      a_k\in\{-1,+1\}~,~~~k=0,1,\dots,n-1~.
$$
The product
$$
    P \= \prod_{k=0}^{n-1} a_k
$$
will be called the parity of the radical.

Obviously, depending on the choice of the $a_k$'s, the pattern may have period less than $n$.
For example, given any $n$,  if $a_k=1$ for all $k$, then the pattern has period 1, giving the well studied radical
{\footnotesize
$$
           \sqrt{2+\sqrt{2+\sqrt{2+\sqrt{2+\cdots}}}}      ~.
$$ 
}
It is easy to find a formula for the number of radicals of minimal period $n$:
Given the radical \eqref{eq:OurRadical}, there are $2^{n}$ different ways to choose the periodic pattern. However,
some of these patterns will have period 1, some period 2, and so on up to period $n$ for all periods $d$ that are divisors of $n$.
Given $n$, let's denote by $N(d)$ the number of radicals with  period $d$. Then
$$
      \sum_{d|n} N(d) \= 2^{n}~.
$$
As is well known, this equation can be inverted with the help of the M\"obius function $\mu(n)$:
$$
      N(n) \= \sum_{d|n}  \mu\lp{n\over d}\rp \, 2^d ~.
$$ 

\section{The Answer.} Towards our result, we present the following lemma \cite{Shklarsky}.
\begin{lemma}
For $\alpha_i\in\{-1,1\}$, $i=0,1,\dots,n-1$,
{\footnotesize
$$
  2\, \sin  \lb \lp \alpha_0+{\alpha_0\alpha_1\over2}+\dots+{\alpha_0\alpha_1\cdots\alpha_{n-1}\over2^{n-1}} \rp {\pi\over4} \rb
  \=
  \alpha_0\sqrt{2+\alpha_1\sqrt{2+\alpha_2\sqrt{2+\dots+\alpha_{n-1}\sqrt{2}}}}~.
$$ 
}
\end{lemma}
The  lemma is easily proven by induction.

According to this lemma, the partial sums of  the continued radical  \eqref{eq:OurRadical} are given by
$$
   x_n \=  2\sin  \lb \lp a_0+{a_0a_1\over2}+\dots+{a_0a_1\cdots a_{n-1}\over2^{n-1}} \rp {\pi\over4} \rb~.
$$
The series
\begin{equation*}
   \alpha_0+{\alpha_0\alpha_1\over2}+\dots+{\alpha_0\alpha_1\cdots\alpha_{n-1}\over2^{n-1}} +\cdots
\end{equation*}
is absolutely convergent and thus it converges to some number $\alpha$. Therefore the original continued radical
 converges to the real number
$$
 x \= 2\sin{\alpha\pi\over4}~.
$$
Alternatively this can be written
$$
 x \= 2\cos{\beta\pi\over2}~,~~~\beta=1-{\alpha\over2}~.
$$
We can find a concise formula for $x$. For this calculation it is more useful to use the products
$$
    P_m \= \prod_{k=0}^{m} a_k~,~~~ m=0,1,\dots,n-1~.
$$
Of course $P_{n-1}=P$ is the parity of the radical. In this notation,
\begin{eqnarray*}
    \alpha  &=&  P_0 + {P_1\over2} + {P_2\over 2^2} + \cdots+  {P_{n-1}\over 2^{n-1}} +   {P_0P_{n-1}\over 2^{n}} +  \cdots  ~,
\end{eqnarray*}
which we can easily rewrite as
\begin{eqnarray*}
   \alpha &=& \lp P_0 + {P_1\over2} + {P_2\over 2^2} + \cdots+  {P_{n-1}\over 2^{n-1}}\rp \, 
                     \lp 1  + {P\over2^n} + {P^2\over 2^{2n}} + \cdots  \rp \\
              &=& \lp P_0 + {P_1\over2} + {P_2\over 2^2} + \cdots +  {P_{n-1}\over 2^{n-1}} \rp \, 
                      {2^n\over 2^n-P}        ~.
\end{eqnarray*}
Then
\begin{eqnarray*}
  {\beta\pi\over2} &=& 2\ell \,  {\pi\over 2^n-P}  ~,
\end{eqnarray*}
where
$$
  2 \ell \=  2^{n-1} -P - (P_{n-2} + 2 P_{n-3} +  \cdots + 2^{n-3}P_1 + 2^{n-2} P_0) ~.
$$

Despite its simplicity, this result may not be easily interpreted. Some examples might help.
In Tables \ref{table:n=2} through \ref{table:n=4} we give the value $x$  of the continued radical for  all possible choices of the $a_k$'s
when $n=2,3$ and 4. The answers take an amazingly compact form.  Looking at these tables and the previous result, we easily realize that  when 
the parity is even then 
$$
     x \= 2 \cos\lp {2\pi \ell \over 2^n -1}\rp~,~~~\ell=0,1,\dots,2^{n-1}-1~,
$$
and when the parity is odd then 
$$
     x \= 2 \cos\lp {2\pi \ell \over 2^n +1}\rp~,~~~\ell=1,\dots,2^{n-1}~.
$$

\begin{table}[ht!]
\begin{center}
\begin{tabular}{|c|c|c|c|c|} \hline
 $\alpha_0$  & $\alpha_1$  & $P$   & ${x\over 2}=\sin{\alpha\pi\over4}$ & ${x\over 2}=\cos{\beta\pi\over2}$ \\ \hline
 $-1$            & $-1$             &+1        &  $\sin\lp-{\pi\over6}\rp$        & $\cos{2\pi\over3}$ \\ \hline
 $-1$            & $+1$            &$-1$   & $\sin\lp-{3\pi\over10}\rp$     & $\cos{4\pi\over5}$ \\ \hline
 $+1$           & $-1$             &$-1$   & $\sin\lp{\pi\over10}\rp$      & $\cos{2\pi\over5}$ \\ \hline
 $+1$           & $+1$            &+1        & $\sin\lp{\pi\over2}\rp$          & $\cos 0$ \\ \hline
\end{tabular}
\end{center}
\caption{\footnotesize The value of our continued radical for all choices of the $a_k$'s when $n=2$. We see that when they 
are chosen such that $\alpha_0\alpha_1=1$, the continued radical equals
$2\cos\lp{2\pi \ell\over 2^2-1}\rp$, $\ell=0,1$, and when they are chosen such that $\alpha_0\alpha_1=-1$, the 
continued radical equals $2\cos\lp{2\pi \ell\over 2^2+1}\rp$, $\ell=1,2$.}
\label{table:n=2}
\end{table}

\begin{table}[hb!]
\begin{center}
\begin{tabular}{|c|c|c|c|c|c|} \hline
 $\alpha_0$  & $\alpha_1$ & $\alpha_2$ &$P$ & ${x\over2}=\sin{\alpha\pi\over4}$ & ${x\over2}=\cos{\beta\pi\over2}$ \\ \hline
 $-1$ & $-1$ & $-1$                                 &$-1$ & $\sin\lp-{\pi\over6}\rp$ & $\cos{6\pi\over9}$ \\ \hline
 $-1$ & $-1$ & $+1$                                &+1  & $\sin\lp-{\pi\over14}\rp$ & $\cos{4\pi\over7}$ \\ \hline
 $-1$ & $+1$ & $-1$                                &+1      & $\sin\lp-{5\pi\over14}\rp$ & $\cos{6\pi\over7}$ \\ \hline
 $-1$ & $+1$ & $+1$                               &$-1$   & $\sin\lp-{7\pi\over18}\rp$ & $\cos{8\pi\over9}$ \\ \hline
 $+1$ & $-1$ & $-1$                                &+1     & $\sin\lp{3\pi\over14}\rp$ & $\cos{2\pi\over7}$ \\ \hline
 $+1$ & $-1$ & $+1$                              &$-1$   & $\sin\lp{\pi\over18}\rp$ & $\cos{4\pi\over9}$ \\ \hline
 $+1$ & $+1$ & $-1$                              &$-1$   & $\sin\lp{5\pi\over18}\rp$ & $\cos{2\pi\over9}$ \\ \hline
 $+1$ & $+1$ & $+1$                             &+1      & $\sin\lp{\pi\over2}\rp$ & $\cos 0$ \\ \hline
\end{tabular}
\end{center}
\caption{\footnotesize The value of our continued radical for all choices of the $a_k$'s when $n=3$. We see that when they 
are chosen such that $\alpha_0\alpha_1\alpha_2=1$, the continued radical equals
$2\cos\lp{2\pi \ell\over 2^3-1}\rp$, $\ell=0,1,2,3$, and when they are chosen such that $\alpha_0\alpha_1\alpha_2=-1$, the 
continued radical equals $2\cos\lp{2\pi \ell\over 2^3+1}\rp$, $\ell=1,2,3,4$.}
\label{table:n=3}
\end{table}

\begin{table}[hb!]
\begin{center}
\begin{tabular}{|c|c|c|c|c|c|c|} \hline
 $\alpha_0$  & $\alpha_1$  & $\alpha_2$ & $\alpha_3$ & $P$ & ${x\over2}=\sin{\alpha\pi\over4}$ & ${x\over2}=\cos{\beta\pi\over2}$ \\ \hline
 $-1$ & $-1$ & $-1$ &  $-1$  & +1 & $\sin\lp-{\pi\over6}\rp$ & $\cos{10\pi\over15}$ \\ \hline
 $-1$ & $-1$ & $-1$ &  +1     &$-1$&$\sin\lp-{7\pi\over34}\rp$ & $\cos{12\pi\over17}$ \\ \hline
 $-1$ & $-1$ & $+1$ & $-1$  & $-1$ &$\sin\lp-{3\pi\over34}\rp$ & $\cos{10\pi\over17}$ \\ \hline
 $-1$ & $-1$ & $+1$ & +1     & +1& $\sin\lp-{\pi\over30}\rp$ & $\cos{8\pi\over15}$ \\ \hline
 $-1$ & $+1$ & $-1$ & $-1$  &$-1$  & $\sin\lp-{11\pi\over34}\rp$ & $\cos{14\pi\over17}$ \\ \hline
 $-1$ & $+1$ & $-1$ & +1     &+1&$\sin\lp-{3\pi\over10}\rp$ & $\cos{12\pi\over15}$ \\ \hline
 $-1$ & $+1$ & $+1$ & $-1$  &+1  & $\sin-\lp{13\pi\over30}\rp$ & $\cos{14\pi\over15}$ \\ \hline
 $-1$ & $+1$ & $+1$ & +1     &$-1$& $\sin\lp-{15\pi\over34}\rp$ & $\cos{16\pi\over17}$ \\ \hline
 $+1$ & $-1$ & $-1$ &  $-1$   &$-1$ & $\sin\lp{5\pi\over34}\rp$ & $\cos{6\pi\over17}$ \\ \hline
 $+1$ & $-1$ & $-1$ &  +1     &+1&$\sin\lp{7\pi\over30}\rp$ & $\cos{4\pi\over15}$ \\ \hline
 $+1$ & $-1$ & $+1$ & $-1$   &+1 &$\sin\lp{\pi\over10}\rp$ & $\cos{6\pi\over15}$ \\ \hline
 $+1$ & $-1$ & $+1$ & +1      &$-1$& $\sin\lp{\pi\over34}\rp$ & $\cos{8\pi\over17}$ \\ \hline
 $+1$ & $+1$ & $-1$ & $-1$   & +1& $\sin\lp{11\pi\over30}\rp$ & $\cos{2\pi\over15}$ \\ \hline
 $+1$ & $+1$ & $-1$ & +1      &$-1$&$\sin\lp{9\pi\over34}\rp$ & $\cos{4\pi\over17}$ \\ \hline
 $+1$ & $+1$ & $+1$ & $-1$   &$-1$ & $\sin\lp{13\pi\over34}\rp$ & $\cos{2\pi\over17}$ \\ \hline
 $+1$ & $+1$ & $+1$ & +1      &+1& $\sin\lp{\pi\over2}\rp$ & $\cos 0$ \\ \hline

\end{tabular}
\end{center}
\caption{\footnotesize The value of our continued radical for all choices of the $a_k$'s when $n=4$. We see that when they 
are chosen such that $\alpha_0\alpha_1\alpha_2\alpha_3=1$, the continued radical equals
$2\cos\lp{2\pi \ell\over 2^4-1}\rp$, $\ell=0,1,2,\dots,7$, and when they are chosen such that $\alpha_0\alpha_1\alpha_2\alpha_3=-1$, the 
continued radical equals $2\cos\lp{2\pi \ell\over 2^4+1}\rp$, $\ell=1,2,\dots,8$.} 
\label{table:n=4}
\end{table}

Motivated by the special cases, we can prove the result in full generality. For $m=0,1,\dots,n-2$ we define 
$$
     Q_m \= {1+P_m\over2}~.
$$
Since $P_m\in\{-1,1\}$, $Q_m\in\{0,1\}$. Inversely, $P_m=2Q_m-1$. We now evaluate $2\ell$ in terms of the $Q_m$'s:
\begin{eqnarray*}
   2\ell &=&  2^{n-1} -P - \lb (2Q_{n-2}-1) + 2(2Q_{n-3}-1) +  \cdots + 2^{n-3}(2Q_1-1) + 2^{n-2}(2Q_0-1) \rb \\
           &=&  2^{n-1} -P + (1+2+\cdots +2^{n-3}+2^{n-2})
                                    -  2(Q_{n-2} + 2Q_{n-3} +  \cdots + 2^{n-3}Q_1 + 2^{n-2}Q_0) \\
           &=&  2^{n-1} -P + (2^{n-1}-1) -2 (Q_{n-2} + 2Q_{n-3} +  \cdots + 2^{n-3}Q_1 + 2^{n-2}Q_0)~,
\end{eqnarray*}
or
$$
      2\ell \=  2^n-P-1-2Q~,
$$
where
$$
     Q \=  Q_{n-2} + 2Q_{n-3} +  \cdots + 2^{n-3}Q_1 + 2^{n-2}Q_0
$$
is the integer whose binary expression is $\overline{Q_0Q_1\cdots Q_{n-3}Q_{n-2}}$. Now we notice that when we go through all possible
sequences $(a_k)_{k=0}^{n-1}$, the sequence $(P_k)_{k=0}^{n-1}$ will go through all possible sequences of $\pm1$'s, and therefore the sequence
$(Q_k)_{k=0}^{n-2}$ will go through all possible sequences of 0's and 1's, with each such sequence appearing once with each value of the parity
$P=P_{n-1}=\pm1$. Consequently, the integer $Q$ will run through the integers from 0 to $2^{n-1}-1$ once with each parity.
Thus, when $P=1$
$$
    \ell \= (2^{n-1}-1)-Q~,
$$
and, as $Q$ runs through the integers from 0 to $2^{n-1}-1$, $\ell$ will run through the same values in reverse order.
When $P=-1$ we get
$$
    \ell \= 2^{n-1}-Q~,
$$
which will then run through all the values from 1 to $2^{n-1}$.

In the following section, we give an alternative way to look at this result:  a nice connection with the Chebycheff polynomials.

\section{Chebycheff Polynomials.}

The $N$th Chebycheff polynomial of the first kind is defined by
$$
    T_{N}(\cos\theta) \= \cos(N\theta) ~.
$$ 
Now consider the quadratic polynomial $P(x)=x^2-2$ defined on $[-2,2]$. Using the substitution
$x=2\cos\theta$, it is easy to see that $P(x)=2\cos(2\theta)$ and
$$
     P^n(x) \= 2\cos(2^n\theta)~.
$$
In other words $P^n(x)=2T_{2^n}(x/2)$. The fixed points of $P^n(x)$ are given by
$P^n(x)=x$, or
$$
     2\cos(2^n\theta) \= 2\cos\theta~.
$$
This equation is easily solved to give the $2^n$ solutions
\begin{eqnarray*}
      \theta &=&  {2\pi \ell\over 2^n-1}~, ~~~\ell\=0,1,\dots,2^{n-1}-1~, \\
      \theta &=&  {2\pi \ell\over 2^n+1}~, ~~~\ell\=1,\dots,2^{n-1}~.
\end{eqnarray*}
The fixed points are then $x=2\cos\theta$.

On the other hand, we can find these fixed points as follows. The equation $P^n(x)=x$ can be written
as $P(P^{n-1}(x))=x$. Using the expression of $P(x)$ we can solve for $P^{n-1}(x)$:
$$
        P^{n-1}(x) \= \pm\sqrt{x+2} ~.
$$
Repeating this $n$ times we find
{\footnotesize
$$
     x\= \pm\sqrt{2\pm\sqrt{2\pm\sqrt{2\pm\cdots\pm\sqrt{2+x}}}}
$$
}
This nested radical reproduces our continued radicals if we iteratively replace  $x$ in the right-hand side by this expression.

\section{Conclusion.}

We have proved that the radicals given by equation \eqref{eq:OurRadical} have limits two times the
fixed points of the Chebycheff polynomials $T_{2^n}(x)$, thus unveiling an interesting relation between these topics.

 In \cite{ZH}, the authors defined the set
$S_2$ of all continued radicals of the form \eqref{eq:OurRadical} (with $a_0=1$) and they investigated some of its properties by assuming that 
the limit of the radicals exists. With this note, we have partially bridged this gap. It is straightforward to see that the limit exists, but we have
identified it only for periodic radicals. 

The continued radical
{\footnotesize
$$
           \sqrt{2+\sqrt{2+\sqrt{2+\sqrt{2+\cdots}}}}     \= 2 
$$ 
}
is well known,  while
{\footnotesize
$$
         \sqrt{2-\sqrt{2+\sqrt{2+\sqrt{2-\cdots}}}}     \= 2 \sin{\pi\over18}
$$ 
}
is a special case of Ramanujan's radical (appearing explicitly in \cite{Ramanujan}).

\section{Acknowledgements.} The author is grateful to the anonymous reviewer for pointing out the articles 
\cite{Nyblom}, \cite{Servi}, and \cite{ZH}.


\bigskip

\noindent\textit{
Department of Physics, University of Central Florida, Orlando, FL 32816 \\
costas@physics.ucf.edu
}


\begin{thebibliography}{99}
    
\bibitem{Andrushkiw}
 R. I. Andrushkiw,
 On the convergence of continued radicals with applications to polynomial equations,
 \textit{ J.  Franklin Inst.} \textbf{319} (1985) 391.
  
  
 \bibitem{Ramanujan}
  B. C. Berndt,
  \textit{Ramanujan's Notebooks}, 
  Springer-Verlag, New York,  1985.
 
  
\bibitem{BB}
  J. M. Borwein and G. de Barra,
  Nested radicals,
  \textit{Amer. Math. Monthly} \textbf{98} (1991) 735.  
  
\bibitem{Dence}  
  T. P. Dence,
  Problem 1174,
  \textit{Math. Mag.} \textbf{56} (1983) 178.
  
\bibitem{Herschfeld}
  A. Herschfeld,
  On infinite radicals,
  \textit{Amer. Math. Monthly} \textbf{42} (1935) 419.
  
  
 \bibitem{JR}
 J. Johnson and T. Richmond,
 Continued radicals,
 \textit{Ramanujan J.} \textbf{15} (2008) 259.
  
   
 \bibitem{Nyblom}
 M. A. Nyblom,
 More nested square roots of 2,
 \textit{Amer. Math. Monthly} \textbf{112} (2005) 822.    
   
   


\bibitem{Shklarsky}
  D. O. Shklarsky, N. N. Chentzov, and I. M. Yaglom,
  \textit{The USSR Problem Book: Selected Problems and Theorems of Elementary Mathematics},
  Dover, New York, 1993.



\bibitem{Servi}
 L. D. Servi, 
 Nested square roots of 2,
 \textit{Amer. Math. Monthly.} \textbf{110} (2003) 326.

\bibitem{Sizer}
  W. S. Sizer,
  Continued roots,
  \textit{Math. Mag.} \textbf{59} (1986) 23.



\bibitem{ZH}
  S. Zimmerman and C. W. Ho,
  On infinitely nested radicals,
  \textit{Math. Mag.} \textbf{81} (2008) 3.  






\end{thebibliography}
\end{document}